\documentclass{qjmam}

\usepackage{graphicx}
\input{amssym.def}
\input{amssym.tex}

\startpage{1}
\yr{}
\vol{}
\issue{}

\def\theequation{\thesection.\arabic{equation}}

\newcommand{\bfa}[1]{\mathbf{#1}}     %
\newcommand{\grad}{\mbox{grad}}       %
\newcommand{\divv}{\mbox{div}}        %
\newcommand{\tgamma}{{\tilde{\gamma}}}%
\newcommand{\jump}[1]{[\![ #1 ]\!]}

\begin{document}
\title[The effect of a curvature-dependent surface tension...]{THE EFFECT OF A CURVATURE-DEPENDENT SURFACE TENSION ON THE SINGULARITIES AT THE TIPS OF A STRAIGHT INTERFACE CRACK}
\author[A.~Y.~Zemlyanova]{A. Y. ZEMLYANOVA}
\address{Department of Mathematics, Texas A\&M University\\
College Station TX 77843}

\received{\recd \revd }

\date{}

\eqnobysec

\maketitle

\noindent

\begin{abstract}
A problem of an interface crack between two semi-planes made out of different materials under an action of an in-plane loading of general tensile-shear type is treated in a semi-analytical manner with the help of Dirichlet-to-Neumann mappings. The boundaries of the crack and the interface between semi-planes are subjected to a curvature-dependent surface tension. The resulting system of six singular integro-differential equations is reduced to the system of three Fredholm equations. It is shown that the introduction of the curvature-dependent surface tension eliminates both classical integrable power singularity of the order $1/2$ and an oscillating singularity present in a classical linear elasticity solutions. The numerical results are obtained by solving the original system of singular integro-differential equations by approximating unknown functions with Taylor polynomials.

\end{abstract}



\setcounter{equation}{0}

\section{Introduction}

The problem of an interface crack is classical for linear elasticity and has generated an enormous amount of literature. The first results in this area have been obtained by Williams \cite{Williams1959} on the example of a semi-infinite crack in dissimilar media. It has been shown for the first time that in addition to the classical singularity of the order $1/2$ the stresses possess an oscillatory singularity of the type $\sin(\epsilon\ln r)$ or $\cos(\epsilon\ln r)$, where $r$ is the radial distance from the tip of the crack. This singularity leads to the physically impossible wrinkling and interpenetration of the banks of the crack near the crack tips. This abnormal behavior has been further confirmed in the studies by  England \cite{England1965}, Erdogan \cite{Erdogan1965}, Malyshev and Salganik \cite{MalyshevSalganik1965},  Rice and Sih \cite{RiceSih1965} and many others. It has also been shown that the zones of overlapping are relatively small ($10^{-7}-10^{-4}$ of the crack length) in the tensile field \cite{Erdogan1965}, however, for the shear loading these zones may be comparable with the size of the crack \cite{ComninouDundurs1980}.

Different approaches have been proposed to overcome the unphysical interpenetration of the crack edges. One of the most well-known models has been proposed by Comninou \cite{Comninou1977a}-\cite{ComninouSchmueser1979}. The main idea is in the introduction of the contact zones near the crack tips. The stress singularity in this case is of the square root type $r^{-1/2}$ and the oscillating singularity disappears. Interface crack problems with contact zones have been also considered in \cite{GautesenDundurs1987}-\cite{Itou1986} and others. Later, Sinclair \cite{Sinclair1980} introduced a crack opening angle in the framework of the Comninou model. 

Another approach has been proposed by Atkinson \cite{Atkinson1977}. In this model the interface is replaced by a thin strip of finite thickness. The crack is either placed inside of the homogeneous strip or is on the boundary of the interface layer and the parameters of the interface layer change continuously. Thus, in both cases the oscillating singularities do not occur. 

Assuming that the interfaces are rough, Mak et al \cite{Mak1980} introduced a no-slip zone near the crack tips. Their solution does not contain oscillating singularities at the crack tips as well. Unlike the fracture model considered by Comninou where a relative shear slip is possible, in this model the relative slip is prohibited and the introduction of the Comninou contact zones is not necessary. In the paper by Boniface and Simha \cite{BonifaceSimha2001} the attempt to suppress an oscillating singularity has been made by choosing a specific crack opening angle which depends on the parameters of interacting dissimilar materials. In another recent study \cite{KimSchiavoneRu2011} the authors eliminated oscillatory singularity by introducing a surface elasticity on the crack boundaries. 

In this paper, we extend the curvature-dependent surface tension approach, first introduced in \cite{OhWaltonSlattery2006}, \cite{SendovaWalton2010} on the example of a straight non-interface crack, to the case of an interface crack between two isotropic semi-planes made out of different materials under an in-plane loading of general type. It is assumed that the equations of linear elasticity are valid everywhere in the bulk in the semi-planes. The classical linear elasticity condition on the crack boundaries and on the interface between semi-planes is augmented with a curvature-dependent surface tension. This problem in a similar formulation has been first approached in \cite{SendovaWalton2010a}. However, this study was incomplete since an essential investigation of the behavior of the stresses and displacements at the crack tips and numerical results were absent. These results are provided in the current paper. Also it is necessary to note, that the regularization of the system of the singular integro-differential equations has been treated in a completely different manner from \cite{SendovaWalton2010a}. The results presented in the current paper show that the introduction of the surface tension on the boundary of an interface crack leads to a complete elimination of both the singularity of the order $1/2$ and the oscillating singularity typical for interface cracks. Thus, the introduction of surface tension removes classical paradoxes of linear elasticity. It is shown that one of the stresses and one of the derivatives of the displacements at the crack tips remains completely bounded, while other parameters assume logarithmic and square-logarithmic singularities. Numerical results presented here are achieved by approximation of the unknown functions by Taylor polynomials. These results are further verified by spline approximations.

\setcounter{equation}{0}

\section{Statement of the problem}

\begin{figure}
	\centering
		\scalebox{0.5}{\includegraphics{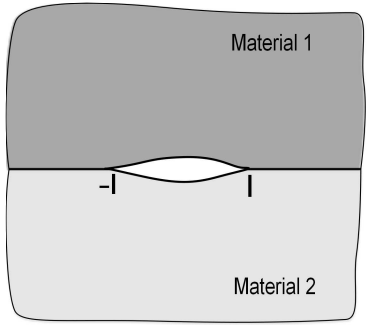}}
	\caption{An interface crack between two dissimilar materials.}
	\label{fig1}
\end{figure}

Consider an upper and a lower semi-planes made out of different materials (fig. \ref{fig1}). The semi-planes are joined perfectly along $x$-axis with the exception of one area of debonding (interface crack) occupying the segment $|x|<l$ along the interface. The materials of the upper and the lower semi-planes have shear moduli $\mu_1$, $\mu_2$ and Poisson's ratios $\nu_1$, $\nu_2$ correspondingly. It is assumed that the constitutive behavior of the material of each of the semi-planes in the bulk can be modeled with the standard equations of linear elasticity and only the boundary conditions are changed to incorporate the curvature dependent surface tension.

Given in-plane stresses $\sigma_{x1}^{\infty}$, $\sigma_{y1}^{\infty}$, $\tau_{xy1}^{\infty}$ and $\sigma_{x2}^{\infty}$, $\sigma_{y2}^{\infty}$, $\tau_{xy2}^{\infty}$ are applied at infinity of each of the semi-planes. The rotations at infinity are $\omega_1^{\infty}$ and $\omega_2^{\infty}$ correspondingly. It is well-known that the stresses at infinity in both planes can not be independent of each other and are connected by the following formulas \cite{RiceSih1965}:
$$
\sigma_{y1}^{\infty}=\sigma_{y2}^{\infty}=\sigma,\,\,\,\,\tau_{xy1}^{\infty}=\tau_{xy2}^{\infty}=\tau,
$$
$$
\frac{1+\kappa_1}{\mu_1}\sigma_{x1}^{\infty}-\frac{1+\kappa_2}{\mu_2}\sigma_{x2}^{\infty}=\left(\frac{3-\kappa_1}{\mu_1}-\frac{3-\kappa_2}{\mu_2}\right)\sigma,
$$
$$
\omega_2^{\infty}-\omega_1^{\infty}=\frac{\mu_2-\mu_1}{2\mu_1\mu_2}\tau.
$$
Here and further the lower index ``1" denotes the parameters belonging to the upper semi-plane and the lower index ``2" denotes correspondingly the parameters of the lower semi-plane. At the upper and the lower boundaries of the crack the tensile $f^{\pm}(x)$ and the shear $g^{\pm}(x)$ stresses are applied. It is assumed that these stresses are in equilibrium:
$$
\int_{-l}^{l}(f^+(x)+ig^+(x))dx-\int_{-l}^l(f^-(x)+ig^-(x))dx=0.
$$

Consider the jump momentum balance condition on the boundaries of the semi-planes \cite{SendovaWalton2010}:
\begin{equation}
\grad_{(\zeta)}\tgamma+2\tgamma H\bfa{n} +\jump{\bfa{T}}\bfa{n} =\bfa{0}.
\label{2_3}
\end{equation}
Here $\tgamma$ is a surface tension on the boundary of the material, $\bfa{T}$ denotes the Cauchy stress tensor, $\bfa{n}$ is the unit normal to the fracture surface $\zeta$ pointing into the bulk of the material, $H=-\frac12\divv_{(\zeta)}\bfa{n}$ is the mean curvature, $\grad_{(\zeta)}$ denotes the surface gradient, and the double brackets $\jump{\ldots}$ denote the jump of the quantity enclosed across the boundary of the semi-plane.

Following the approach proposed in \cite{SendovaWalton2010} for a straight non-interface crack, assume that the curvature dependent surface tension acts on the crack boundaries and on the interface between two different materials and has the following form:
\begin{equation}
\tilde{\gamma}^{\pm}=\gamma_0^{\pm}+\gamma_1^{\pm}\divv_{(\zeta)}{\bfa{n}},\,\,\,\,|x|<l,
\label{2_1}
\end{equation}
\begin{equation}
\tilde{\gamma}^{i}=\gamma_0^{i}+\gamma_1^{i}\divv_{(\zeta)}{\bfa{n}},\,\,\,\,|x|>l,
\label{2_2}
\end{equation}
where $\gamma_0^{\pm}$, $\gamma_1^{\pm}$, $\gamma_0^i$, $\gamma_1^i$ are real constants, the upper indexes ``$+$", ``$-$" and ``$i$" describe the parameters on the upper bank of the crack, on the lower bank of the crack and on the interface between the materials outside of the crack correspondingly.

Assuming that the derivatives of the displacements $u_{1,1}$, $u_{2,1}$ and their derivatives are small, compute the surface divergence $\divv_{(\zeta)}$ of the inner normal $\bfa{n}$ \cite{Slatteryetal}:
\begin{equation}
\divv_{(\zeta)}\bfa{n}=u_{2,11}(x).
\label{2_4}
\end{equation}

Substituting the conditions (\ref{2_1}), (\ref{2_2}) into the condition (\ref{2_3}) and taking into account (\ref{2_4}), one obtains the following boundary conditions:
\begin{equation}
\sigma_{12}^+(x,0)=\gamma_1^+u_{2,111}^+(x,0)+f^+(x),\,\,\,|x|<l,
\label{2_5}
\end{equation}
$$
\sigma_{22}^+(x,0)=-\gamma_0^+u_{2,11}^+(x,0)+g^+(x),\,\,\,|x|<l,
$$
on the upper boundary of the crack;
\begin{equation}
\sigma_{12}^-(x,0)=-\gamma_1^-u_{2,111}^-(x,0)+f^-(x),\,\,\,|x|<l,
\label{2_6}
\end{equation}
$$
\sigma_{22}^-(x,0)=\gamma_0^-u_{2,11}^-(x,0)+g^-(x),\,\,\,|x|<l,
$$
on the lower boundary of the crack; and
\begin{equation}
\sigma_{12}^+(x,0)-\sigma_{12}^-(x,0)=\gamma_1^{i}u_{2,111}(x,0),\,\,\,|x|>l,
\label{2_7}
\end{equation}
$$
\sigma_{22}^+(x,0)-\sigma_{22}^-(x,0)=-\gamma_0^{i}u_{2,11}(x,0),\,\,\,|x|>l,
$$
on the interface between two semi-planes. Here and further ``$+$" and ``$-$" signs describe the limiting values of the stresses and the displacements on the boundaries of the semi-planes from the side of the upper and the lower semi-planes.  

Observe that due to the perfect attachment of the semi-planes along the interface the displacements $u_1$, $u_2$, and hence their derivatives, are equal on the interface. Thus, these functions can be taken without ``$\pm$" signs:
$$
u_{1,1}^+(x,0)=u_{1,1}^-(x,0)=u_{1,1}(x,0),\,\,\,u_{2,1}^+(x,0)=u_{2,1}^-(x,0)=u_{2,1}(x,0),\,\,\,|x|>l.
$$

\setcounter{equation}{0}

\section{Reduction to the system of integro-differential equations}

The Dirichlet-to-Neumann transforms will be used to solve the problem. It has been shown \cite{SendovaWalton2010} that the Dirichlet-to-Neumann mappings for the upper semi-plane have a form:
\begin{equation}
\sigma_{12}(x,0)=\alpha_1u_{2,1}(x,0)+\beta_1\int_{-\infty}^{\infty}\frac{u_{1,1}(r,0)dr}{r-x};
\label{3_1}
\end{equation}
$$
\sigma_{22}(x,0)=-\alpha_1u_{1,1}(x,0)+\beta_1\int_{-\infty}^{\infty}\frac{u_{2,1}(r,0)dr}{r-x}.
$$
In a similar way as in \cite{SendovaWalton2010} it is possible to derive that the Dirichlet-to-Neumann mappings for the lower semi-plane are given by the formulas:
\begin{equation}
\sigma_{12}(x,0)=\alpha_2u_{2,1}(x,0)-\beta_2\int_{-\infty}^{\infty}\frac{u_{1,1}(r,0)dr}{r-x};
\label{3_2}
\end{equation}
$$
\sigma_{22}(x,0)=-\alpha_2u_{1,1}(x,0)-\beta_2\int_{-\infty}^{\infty}\frac{u_{2,1}(r,0)dr}{r-x},
$$
where $\alpha_1$, $\alpha_2$, $\beta_1$ and $\beta_2$ are given constants which depend on the material parameters of the semi-planes:
$$
\alpha_1=\frac{2\mu_1^2}{\lambda_1+3\mu_1},\,\,\,\alpha_2=\frac{2\mu_2^2}{\lambda_2+3\mu_2},\,\,\,\beta_1=\frac{2\mu_1(\lambda_1+2\mu_1)}{(\lambda_1+3\mu_1)\pi},\,\,\,\beta_2=\frac{2\mu_2(\lambda_2+2\mu_2)}{(\lambda_2+3\mu_2)\pi},
$$
and $\nu_1=\frac{\lambda_1}{2(\lambda_1+\mu_1)}$, $\nu_2=\frac{\lambda_2}{2(\lambda_2+\mu_2)}$.

Since the formulas (\ref{3_1}), (\ref{3_2}) are derived under the assumption that the far field loading in the semi-planes vanishes at infinity, it is necessary to separate the homogeneous distant loading field generated by the stresses applied at infinity. This leads to the following modifications to the boundary conditions (\ref{2_5}), (\ref{2_6}) and (\ref{2_7}): 
\begin{equation}
\sigma_{12}^{0+}(x,0)+\tau=\gamma_1^+u_{2,111}^{0+}(x,0)+f^+(x),\,\,\,|x|<l,
\label{3_3}
\end{equation}
$$
\sigma_{22}^{0+}(x,0)+\sigma=-\gamma_0^+u_{2,11}^{0+}(x,0)+g^+(x),\,\,\,|x|<l,
$$
\begin{equation}
\sigma_{12}^{0-}(x,0)+\tau=-\gamma_1^-u_{2,111}^{0-}(x,0)+f^-(x),\,\,\,|x|<l,
\label{3_4}
\end{equation}
$$
\sigma_{22}^{0-}(x,0)+\sigma=\gamma_0^-u_{2,11}^{0-}(x,0)+g^-(x),\,\,\,|x|<l.
$$
\begin{equation}
\sigma_{12}^{0+}(x,0)-\sigma_{12}^{0-}(x,0)=\gamma_1^{i}u_{2,111}^0(x,0),\,\,\,|x|>l,
\label{3_5}
\end{equation}
$$
\sigma_{22}^{0+}(x,0)-\sigma_{22}^{0-}(x,0)=-\gamma_0^{i}u_{2,11}^0(x,0),\,\,\,|x|>l,
$$
where all the parameters with the upper index ``$0$" denote the stresses and the displacements in the absence of the loading at infinity.

Introduce new unknown functions:
\begin{equation}
\varphi_1(x)=u_{1,1}^{0+}(x,0),\,\,\,\varphi_2(x)=u_{1,1}^{0-}(x,0),\,\,\,|x|<l,\,\,\,\varphi(x)=u_{1,1}^0(x,0),\,\,\,|x|>l,
\label{3_6}
\end{equation}
$$
\psi_1(x)=u_{2,1}^{0+}(x,0),\,\,\,\psi_2(x)=u_{2,1}^{0-}(x,0),\,\,\,|x|<l,\,\,\,\psi(x)=u_{2,1}^0(x,0),\,\,\,|x|>l.
$$

Substituting the Dirichlet-to-Neumann mappings (\ref{3_1}), (\ref{3_2}) into the boundary conditions (\ref{3_3}), (\ref{3_4}) and (\ref{3_5}) leads to the following system of six singular integro-differential equations for six unknown functions (\ref{3_6}):
$$
\alpha_1\psi_1(x)+\beta_1\int_{|r|<l}\frac{\varphi_1(r)dr}{r-x}+\beta_1\int_{|r|>l}\frac{\varphi(r)dr}{r-x}+\tau=\gamma_1^+\psi_1''(x)+f^+(x),\,\,\,|x|<l,
$$
$$
-\alpha_1\varphi_1(x)+\beta_1\int_{|r|<l}\frac{\psi_1(r)dr}{r-x}+\beta_1\int_{|r|>l}\frac{\psi(r)dr}{r-x}+\sigma=-\gamma_0^+\psi_1'(x)+g^+(x),\,\,\,|x|<l,
$$
$$
\alpha_2\psi_2(x)-\beta_2\int_{|r|<l}\frac{\varphi_2(r)dr}{r-x}-\beta_2\int_{|r|>l}\frac{\varphi(r)dr}{r-x}+\tau=-\gamma_1^-\psi_2''(x)+f^-(x),\,\,\,|x|<l,
$$
\begin{equation}
-\alpha_2\varphi_2(x)-\beta_2\int_{|r|<l}\frac{\psi_2(r)dr}{r-x}-\beta_2\int_{|r|>l}\frac{\psi(r)dr}{r-x}+\sigma=\gamma_0^-\psi_2'(x)+g^-(x),\,\,\,|x|<l,
\label{3_7}
\end{equation}
$$
(\alpha_1-\alpha_2)\psi(x)+\beta_1\int_{|r|<l}\frac{\varphi_1(r)dr}{r-x}+\beta_2\int_{|r|<l}\frac{\varphi_2(r)dr}{r-x}+(\beta_1+\beta_2)\int_{|r|>l}\frac{\varphi(r)dr}{r-x}
$$
$$
=\gamma_1^i\psi''(x),\,\,\,|x|>l,
$$
$$
-(\alpha_1-\alpha_2)\varphi(x)+\beta_1\int_{|r|<l}\frac{\psi_1(r)dr}{r-x}+\beta_2\int_{|r|<l}\frac{\psi_2(r)dr}{r-x}+(\beta_1+\beta_2)\int_{|r|>l}\frac{\psi(r)dr}{r-x}
$$
$$
=-\gamma_0^i\psi'(x),\,\,\,|x|>l.
$$

Observe that the unknown functions $\varphi_1(x)$, $\varphi_2(x)$ and $\varphi(x)$ can be easily eliminated from the system of the integral equations (\ref{3_7}). To achieve this integrate the second, the fourth and the sixth equations of the system (\ref{3_7}) with the Cauchy kernel and change the order of integration using the Poincare-Bertrand formula \cite{Gakhov1990}:
\begin{equation}
-\alpha_1\int_{|r|<l}\frac{\varphi_1(r)dr}{r-x}-\beta_1\pi^2\psi_1(x)+\beta_1\int_{|r|<l}\ln\left|\frac{(l-r)(l+x)}{(l-x)(l+r)} \right|\frac{\psi_1(r)dr}{r-x}
\label{3_8}
\end{equation}
$$
+\beta_1\int_{|r|>l}\ln\left|\frac{(l-r)(l+x)}{(l-x)(l+r)} \right|\frac{\psi(r)dr}{r-x}+\gamma_0^+\int_{|r|<l}\frac{\psi'_1(r)dr}{r-x}=-\sigma\ln\left|\frac{l-x}{l+x} \right|
$$
$$
+\int_{|r|<l}\frac{g^+(r)dr}{r-x},\,\,\,|x|<l, 
$$
\begin{equation}
-\alpha_2\int_{|r|<l}\frac{\varphi_2(r)dr}{r-x}+\beta_2\pi^2\psi_2(x)-\beta_2\int_{|r|<l}\ln\left|\frac{(l-r)(l+x)}{(l-x)(l+r)} \right|\frac{\psi_2(r)dr}{r-x}
\label{3_9}
\end{equation}
$$
-\beta_2\int_{|r|>l}\ln\left|\frac{(l-r)(l+x)}{(l-x)(l+r)} \right|\frac{\psi(r)dr}{r-x}-\gamma_0^-\int_{|r|<l}\frac{\psi'_2(r)dr}{r-x}=-\sigma\ln\left|\frac{l-x}{l+x} \right|
$$
$$
+\int_{|r|<l}\frac{g^-(r)dr}{r-x},\,\,\,|x|<l, 
$$
\begin{equation}
-(\alpha_1-\alpha_2)\int_{|r|>l}\frac{\varphi(r)dr}{r-x}-(\beta_1+\beta_2)\pi^2\psi(x)-\beta_1\int_{|r|<l}\ln\left|\frac{(l-r)(l+x)}{(l-x)(l+r)} \right|\frac{\psi_1(r)dr}{r-x}
\label{3_10}
\end{equation}
$$
-\beta_2\int_{|r|<l}\ln\left|\frac{(l-r)(l+x)}{(l-x)(l+r)} \right|\frac{\psi_2(r)dr}{r-x}-(\beta_1+\beta_2)\int_{|r|>l}\ln\left|\frac{(l-r)(l+x)}{(l-x)(l+r)} \right|\frac{\psi(r)dr}{r-x}
$$
$$
+\gamma_0^i\int_{|r|<l}\frac{\psi'(r)dr}{r-x}=0,\,\,\,|x|>l. 
$$
Observe that for $|x|>l$ the formulas (\ref{3_8}) and (\ref{3_9}) look exactly the same with the exception that the terms $-\beta_1\pi^2\psi_1(x)$ and $\beta_2\pi^2\psi_2(x)$ drop out. Similarly, for $|x|<l$ the formula (\ref{3_10}) is lacking the term $-(\beta_1+\beta_2)\pi^2\psi(x)$. These are the consequences of the fact that the order of integrals can be interchanged if at least one of the kernels is regular.

Substituting the formulas (\ref{3_8}), (\ref{3_9}) and (\ref{3_10}) into the first, the third and the fifth equations of the system (\ref{3_7}) leads to the following system of three integro-differential equations:
$$
-\gamma_1^+\psi''_1(x)+\alpha_1\left(1-\frac{\beta_1^2\pi^2}{\alpha_1^2} \right)\psi_1(x)+\frac{\gamma_0^+\beta_1}{\alpha_1}\int_{|r|<l}\frac{\psi'_1(r)dr}{r-x}+\frac{\gamma_0^i\beta_1}{\alpha_1-\alpha_2}\int_{|r|>l}\frac{\psi'(r)dr}{r-x}
$$
$$
-\frac{\beta_1^2\alpha_2}{\alpha_1(\alpha_1-\alpha_2)}\int_{|r|<l}\ln\left|\frac{(l-r)(l+x)}{(l-x)(l+r)} \right|\frac{\psi_1(r)dr}{r-x}-\frac{\beta_1\beta_2}{\alpha_1-\alpha_2}\int_{|r|<l}\ln\left|\frac{(l-r)(l+x)}{(l-x)(l+r)} \right|\frac{\psi_2(r)dr}{r-x}
$$
$$
-\frac{\beta_1(\beta_1\alpha_2+\beta_2\alpha_1)}{\alpha_1(\alpha_1-\alpha_2)}\int_{|r|>l}\ln\left|\frac{(l-r)(l+x)}{(l-x)(l+r)} \right|\frac{\psi(r)dr}{r-x}
$$
$$
=\frac{\beta_1}{\alpha_1}\int_{|r|<l}\frac{g^+(r)dr}{r-x}+f^+(x)-\sigma\frac{\beta_1}{\alpha_1}\ln\left|\frac{l-x}{l+x} \right|-\tau, \,\,\, |x|<l,
$$
$$
\gamma_1^-\psi''_2(x)+\alpha_2\left(1-\frac{\beta_2^2\pi^2}{\alpha_2^2} \right)\psi_2(x)+\frac{\gamma_0^-\beta_2}{\alpha_2}\int_{|r|<l}\frac{\psi'_2(r)dr}{r-x}-\frac{\gamma_0^i\beta_2}{\alpha_1-\alpha_2}\int_{|r|>l}\frac{\psi'(r)dr}{r-x}
$$
$$
+\frac{\beta_1\beta_2}{\alpha_1-\alpha_2}\int_{|r|<l}\ln\left|\frac{(l-r)(l+x)}{(l-x)(l+r)} \right|\frac{\psi_1(r)dr}{r-x}+\frac{\beta_2^2\alpha_1}{\alpha_2(\alpha_1-\alpha_2)}\int_{|r|<l}\ln\left|\frac{(l-r)(l+x)}{(l-x)(l+r)} \right|\frac{\psi_2(r)dr}{r-x}
$$
\begin{equation}
+\frac{\beta_2(\beta_1\alpha_2+\beta_2\alpha_1)}{\alpha_2(\alpha_1-\alpha_2)}\int_{|r|>l}\ln\left|\frac{(l-r)(l+x)}{(l-x)(l+r)} \right|\frac{\psi(r)dr}{r-x}
\label{3_11}
\end{equation}
$$
=-\frac{\beta_2}{\alpha_2}\int_{|r|<l}\frac{g^-(r)dr}{r-x}+f^-(x)+\sigma\frac{\beta_2}{\alpha_2}\ln\left|\frac{l-x}{l+x} \right|-\tau, \,\,\, |x|<l,
$$
$$
-\gamma_1^i\psi''(x)+(\alpha_1-\alpha_2)\left(1-\frac{(\beta_1+\beta_2)^2\pi^2}{(\alpha_1-\alpha_2)^2} \right)\psi(x)+\frac{\gamma_0^+\beta_1}{\alpha_1}\int_{|r|<l}\frac{\psi'_1(r)dr}{r-x}
$$
$$
-\frac{\gamma_0^-\beta_2}{\alpha_2}\int_{|r|<l}\frac{\psi'_2(r)dr}{r-x}+\frac{\gamma_0^i(\beta_1+\beta_2)}{\alpha_1-\alpha_2}\int_{|r|>l}\frac{\psi'(r)dr}{r-x}
$$
$$
-\frac{\beta_1(\beta_1\alpha_2+\beta_2\alpha_1)}{\alpha_1(\alpha_1-\alpha_2)}\int_{|r|<l}\ln\left|\frac{(l-r)(l+x)}{(l-x)(l+r)}\right|\frac{\psi_1(r)dr}{r-x}
$$
$$
-\frac{\beta_2(\beta_1\alpha_2+\beta_2\alpha_1)}{\alpha_2(\alpha_1-\alpha_2)}\int_{|r|<l}\ln\left|\frac{(l-r)(l+x)}{(l-x)(l+r)} \right|\frac{\psi_2(r)dr}{r-x}
$$
$$
-\frac{(\beta_1\alpha_2+\beta_2\alpha_1)^2}{\alpha_1\alpha_2(\alpha_1-\alpha_2)}\int_{|r|>l}\ln\left|\frac{(l-r)(l+x)}{(l-x)(l+r)} \right|\frac{\psi(r)dr}{r-x}
$$
$$
=\frac{\beta_1}{\alpha_1}\int_{|r|<l}\frac{g^+(r)dr}{r-x}+\frac{\beta_2}{\alpha_2}\int_{|r|<l}\frac{g^-(r)dr}{r-x}-\sigma\left(\frac{\beta_1}{\alpha_1}+\frac{\beta_2}{\alpha_2}\right)\ln\left|\frac{l-x}{l+x} \right|, \,\,\, |x|>l,
$$
assuming that $\alpha_1\ne \alpha_2$.

Observe that the system (\ref{3_11}) is not a system of singular integro-differential equations but rather the system of integro-differential equations with logarithmic singularities.

Finally, to completely regularize the system in the spirit of \cite{MikhPros1986}, introduce the following new unknown functions:
\begin{equation}
\chi_1(x)=\psi''_1(x),\,\,\,\,\,\chi_2(x)=\psi''_2(x),\,\,\,\,\,\chi(x)=\psi''(x).
\label{3_12}
\end{equation}
By integrating the formulas (\ref{3_12}), one obtains:
$$
\psi'_1(x)=\int_{|r|<l}\omega_0(r,x)\chi_1(r)dr+C_1,
$$
$$
\psi_1(x)=\int_{|r|<l}\omega_1(r,x)\chi_1(r)dr+C_1x+C_2,
$$
$$
\psi'_2(x)=\int_{|r|<l}\omega_0(r,x)\chi_2(r)dr+C_3,
$$
$$
\psi_2(x)=\int_{|r|<l}\omega_1(r,x)\chi_2(r)dr+C_3x+C_4,
$$
\begin{equation}
\psi'(x)=\int_{-\infty}^x\chi(r)dr,\,\,\,x<-l,
\label{3_13}
\end{equation}
$$
\psi(x)=\int_{-\infty}^x(x-r)\chi(r)dr,\,\,\,x<-l,
$$
$$
\psi'(x)=-\int_x^{\infty}\chi(r)dr,\,\,\,x>l,
$$
$$
\psi(x)=-\int_x^{\infty}(x-r)\chi(r)dr,\,\,\,x>l,
$$
where
$$
\omega_0(r,x)=\left\{
\begin{array}{ll}
1, & r \in [-l,x],\\
0, & r \notin [-l,x],
\end{array}
\right.\,\,\,\,\,
\omega_1(r,x)=(x-r)\omega_0(r,x).
$$

The representations (\ref{3_13}) contain four real constants $C_1$, $C_2$, $C_3$ and $C_4$. These constants will be found from additional physical conditions which will be stated later. The corresponding integration constants for the function $\chi(x)$ are chosen to be equal to zero to guarantee the appropriate decay of the function $\psi(x)$ at infinity and the existence of all the integrals in (\ref{3_11}).

Finally, substituting the formulas (\ref{3_12}), (\ref{3_13}) into the system (\ref{3_11}) reduces this system to the final system of Fredholm equations the theory of which is well-developed.

\setcounter{equation}{0}

\section{Additional conditions}

The representations (\ref{3_13}) contain four real constants $C_1$, $C_2$, $C_3$ and $C_4$. To fix these constants four additional real conditions coming from the physics of the problem will be stated. Firstly, it is reasonable to assume that the total stresses applied to the boundaries of the crack must be equal to zero, which gives two real conditions:
\begin{equation}
\int_{-l}^l (\sigma_{22}^+-\sigma_{22}^-)dx=0,
\label{4_1}
\end{equation}
\begin{equation}
\int_{-l}^l (\sigma_{12}^+-\sigma_{12}^-)dx=0.
\label{4_2}
\end{equation}
The equation (\ref{4_1}) can be expressed through the unknown functions and the constants $C_j$ using the boundary conditions (\ref{3_3}), (\ref{3_4}):
$$
\int_{-l}^l(\gamma_0^+\psi'_1(x)+\gamma_0^-\psi'_2(x)-g^+(x)+g^-(x))dx=0.
$$
At the same time to rewrite the condition (\ref{4_2}) in terms of the unknown functions and the constants $C_j$ it is better to use the Dirichlet-to-Neumann mappings (\ref{3_1}), (\ref{3_2}). 

To make the solution physically reasonable it is necessary for the displacements to be single-valued along the boundary of the crack. In other words, if the position of the left end of the crack is fixed then the displacement of the right end of the crack must be the same along both the upper and the lower semi-planes. This leads to the following two real conditions:
\begin{equation}
\int_{-l}^l (u^+_{1,1}-u^-_{1,1})dx=0,
\label{4_3}
\end{equation}
\begin{equation}
\int_{-l}^l (u^+_{2,1}-u^-_{2,1})dx=0.
\label{4_4}
\end{equation}
These conditions are easy to rewrite in terms of the unknown functions and the constants $C_j$ as:
$$
\int_{-l}^l(\varphi_1(x)-\varphi_2(x))dx=0,
$$
$$
\int_{-l}^l(\psi_1(x)-\psi_2(x))dx=0.
$$

The conditions (\ref{4_1})-(\ref{4_4}) allow us to fix all the constants in the formulas (\ref{3_13}).

\setcounter{equation}{0}

\section{Singularities of the stresses and derivatives of the displacements at the crack tips}

Consider the system of singular integro-differential equations (\ref{3_9}). Observe that the functions $\varphi_1(x)$, $\varphi_2(x)$, $\varphi(x)$, $\psi_1(x)$, $\psi_2(x)$ and $\psi(x)$ can not have integrable power singularities of any power (including pure imaginary powers) at the crack tips $x=\pm l$. Indeed, assume that the singularities are possible and consider temporarily only the point $x=l$. Then 
\begin{equation}
\varphi_k(x)=A_k(x-l)^{\gamma}+o((x-l)^{\gamma}),
\label{5_1}
\end{equation}
$$
\psi''_k(x)=B_k(x-l)^{\gamma}+o((x-l)^{\gamma}),
$$
where the index $k$ is equal to $1$, $2$ or missing and $A_k$, $B_k$ are real coefficients. It follows then that the functions $\psi'_k(x)$, $\psi_k(x)$ are bounded at the point $x=l$. From the second, the fourth and the sixth equations of the system (\ref{3_9}) it follows then that the functions $\varphi_k(x)$ may have at most logarithmic singularities. Using this fact and the first, the third and the fifth equations of the system (\ref{3_9}), deduce then that the functions $\psi''_k(x)$ may have at most logarithmic or square-logarithmic singularities. Thus, we have a contradiction with the initial representations (\ref{5_1}).

Hence, the functions $\varphi_k(x)$, $\psi''_k(x)$ can not have integrable power singularities of any power. This means, that unlike in the linear theory of elastic fracture mechanics (LEFM), the stresses and the derivatives of the displacements do not possess singularities of the order $1/2$ or oscillating singularities with a pure imaginary power. Thus, the introduction of the surface tension on the boundary of the crack removes these important contradictions of the linear theory.

On the other hand, it can be seen from the system (\ref{3_9}) that the functions $\varphi_k(x)$, $\psi''_k(x)$ can have logarithmic singularities of the following form:
\begin{equation}
\varphi_k(x)=A_k\ln(x-l)+O(1),
\label{5_2}
\end{equation}
$$
\psi''_k(x)=B_k\ln^2(x-l)+C_k\ln(x-l)+O(1),
$$
where $A_k$, $B_k$ and $C_k$ are real coefficients. Then, as before it follows that the functions $\psi'_k(x)$, $\psi_k(x)$ are bounded at the point $x=l$. Hence, using the boundary conditions (\ref{3_3}), (\ref{3_4}), (\ref{3_5}) it can be seen that the stresses $\sigma_{22}^{\pm}$ are bounded at the point $x=l$ and the stresses $\sigma_{12}^{\pm}$ have a logarithmic singularity of the same type as the functions $\psi''_k(x)$ (\ref{5_2}). The derivative of the displacement $u_{2,1}^{\pm}$ is bounded and the derivative $u_{1,1}^{\pm}$ has the singularity of the same type as $\varphi_k(x)$. Obviously, the same conclusions are valid at the left tip of the crack $x=-l$. Thus, the singularities of the stresses and of the displacements are of the same type as for a curvilinear non-interface crack \cite{ZemWal} and hence, the surface tension model does not differentiate between an interface and a non-interface cracks. Observe also, that for the conclusion of this section to stand it is essential that the coefficients $\gamma_1^{\pm}$ and $\gamma_1^i$ are all non-zero.

\setcounter{equation}{0}

\section{Numerical results}

To produce the numerical results for the stresses and the displacements, the singular integro-differential system (\ref{3_7}) which has a significantly more convenient form than the regularized system (\ref{3_11}) will be solved numerically. Two different approaches will be considered. The first approach consists in representing the unknown functions by Taylor polynomials, substituting these polynomials into the system (\ref{3_7}) and consequently solving the resulting system of linear algebraic equations. This approach has been found to be the most convenient and to produce the best numerical results. For a verification purposes only this approach has been corroborated by representing the unknown functions with splines.

First of all, the unknown functions $\varphi_k(x)$, $\psi_k(x)$ are approximated by Taylor polynomials with unknown coefficients:
\begin{equation}
\varphi_1(x)=\sum_{k=0}^{N}a_k^1x^k,\,\,\,\varphi_2(x)=\sum_{k=0}^{N}a_k^2x^k,\,\,\,\varphi(x)=\sum_{k=1}^{N+1}a_k^3x^{-k},
\label{6_1}
\end{equation}
$$
\psi_1(x)=\sum_{k=0}^{N+2}b_k^1x^k,\,\,\,\psi_2(x)=\sum_{k=0}^{N+2}b_k^2x^k,\,\,\,\psi(x)=\sum_{k=1}^{N+1}b_k^3x^{-k}.
$$
The idea behind the representations (\ref{6_1}) is that all of the series have the same number of terms ($N+1$) except for the series for $\psi_1(x)$ and $\psi_2(x)$ which have two additional terms which correspond to the constants $C_1$, $C_2$, $C_3$ and $C_4$. In addition, the functions $\varphi(x)$, $\psi(x)$ must decline to zero at infinity which guarantees the convergence of all the integrals involved.

Substituting the representations (\ref{6_1}) into the system (\ref{3_7}) leads to the following integrals which are easy to compute analytically:
$$
\int_{|r|<l}\frac{r^kdr}{r-x}=-\sum_{j=0,j\neq k}^{\infty}\frac{1-(-1)^{j-k}}{j-k}x^jl^{k-j},\,\,\,|x|<l,
$$
\begin{equation}
\int_{|r|>l}\frac{r^{-k}dr}{r-x}=\sum_{j=0}^{\infty}\frac{1-(-1)^{j+k}}{j+k}x^jl^{-j-k},\,\,\,|x|<l,
\label{6_2}
\end{equation}
$$
\int_{|r|<l}\frac{r^kdr}{r-x}=-\sum_{j=1}^{\infty}\frac{1-(-1)^{j+k}}{j+k}x^{-j}l^{j+k},\,\,\,|x|>l,
$$
$$
\int_{|r|>l}\frac{r^{-k}dr}{r-x}=\sum_{j=1,j\neq k}^{\infty}\frac{1-(-1)^{j-k}}{j-k}x^{-j}l^{j-k},\,\,\,|x|>l.
$$

Assume that the tensile and the shear stresses applied to the boundaries of the crack have the following Taylor series representations:
\begin{equation}
f^+(x)=\sum_{k=0}^Nf_k^+x^k,\,\,\,f^-(x)=\sum_{k=0}^Nf_k^-x^k,
\label{6_3}
\end{equation}
$$
g^+(x)=\sum_{k=0}^Ng_k^+x^k,\,\,\,g^-(x)=\sum_{k=0}^Ng_k^-x^k.
$$

By substituting the formulas (\ref{6_1}), (\ref{6_2}), (\ref{6_3}) into the system (\ref{3_7}) and making sure that the coefficients by the powers of $x$ are the same on both sides of the equations, obtain the following system of linear algebraic equations:
$$
\alpha_1b_k^1+\beta_1\left(-\sum_{j=0,j\neq k}^{N}\frac{1-(-1)^{k-j}}{k-j}a_j^1l^{j-k}+\sum_{j=1}^{N+1}\frac{1-(-1)^{j+k}}{j+k}a_j^3l^{-k-j} \right)
$$
$$
-\gamma_1^+(k+2)(k+1)b_{k+2}^1=-\tau\delta_{0k}+f_k^+,\,\,\,k=0,\ldots,N,
$$
$$
-\alpha_1a_k^1+\beta_1\left(-\sum_{j=0,j\neq k}^{N+2}\frac{1-(-1)^{k-j}}{k-j}b_j^1l^{j-k}+\sum_{j=1}^{N+1}\frac{1-(-1)^{j+k}}{j+k}b_j^3l^{-k-j} \right)
$$
$$
+\gamma_0^+(k+1)b_{k+1}^1=-\sigma\delta_{0k}+g_k^+,\,\,\,k=0,\ldots,N,
$$
$$
\alpha_2b_k^2-\beta_2\left(-\sum_{j=0,j\neq k}^{N}\frac{1-(-1)^{k-j}}{k-j}a_j^2l^{j-k}+\sum_{j=1}^{N+1}\frac{1-(-1)^{j+k}}{j+k}b_j^3l^{-k-j} \right)
$$
$$
+\gamma_1^-(k+2)(k+1)b_{k+2}^2=-\tau\delta_{0k}+f_k^-,\,\,\,k=0,\ldots,N,
$$
\begin{equation}
-\alpha_2a_k^2-\beta_2\left(-\sum_{j=0,j\neq k}^{N+2}\frac{1-(-1)^{k-j}}{k-j}b_j^2l^{j-k}+\sum_{j=1}^{N+1}\frac{1-(-1)^{j+k}}{j+k}b_j^3l^{-k-j} \right)
\label{6_4}
\end{equation}
$$
-\gamma_0^-(k+1)b_{k+1}^2=-\sigma\delta_{0k}+g_k^-,\,\,\,k=0,\ldots,N,
$$
$$
(\alpha_1-\alpha_2)b_k^3+(\beta_1+\beta_2)\sum_{j=1,j\neq k}^{N+1}\frac{1-(-1)^{k-j}}{k-j}a_j^3l^{k-j}-\beta_1\sum_{j=0}^{N}\frac{1-(-1)^{j+k}}{j+k}b_j^1l^{k+j}
$$
$$
-\beta_2\sum_{j=0}^{N}\frac{1-(-1)^{j+k}}{j+k}b_j^2l^{k+j}-\gamma_1^i(k-2)(k-1)b_{k-2}^3=0,\,\,\,k=1,\ldots,N+1,
$$
$$
-(\alpha_1-\alpha_2)a_k^3+(\beta_1+\beta_2)\sum_{j=1,j\neq k}^{N+1}\frac{1-(-1)^{k-j}}{k-j}b_j^3l^{k-j}-\beta_1\sum_{j=0}^{N+2}\frac{1-(-1)^{j+k}}{j+k}b_j^1l^{j+k} 
$$
$$
-\beta_2\sum_{j=0}^{N+2}\frac{1-(-1)^{j+k}}{j+k}b_j^2l^{j+k}-\gamma_0^i(k-1)b_{k-1}^3=0,\,\,\,k=1,\ldots,N+1,
$$
where $\delta_{jk}$ is the Kronecker symbol.

It is also possible to rewrite the additional conditions (\ref{4_1})-(\ref{4_4}) in terms of the Taylor series representations (\ref{6_1}). Thus, the following four additional equations for the coefficients $a_j^k$, $b_j^k$ must be satisfied:
\begin{equation}
\gamma_0^+\sum_{k=0}^{N+2}b_k^1e^k(1-(-1)^k)+\gamma_0^-\sum_{k=0}^{N+2}b_k^2e^k(1-(-1)^k)=0,
\label{6_5}
\end{equation}
\begin{equation}
\gamma_1^+\sum_{k=1}^{N+2}b_k^1ke^{k-1}(1-(-1)^{k-1})+\gamma_1^-\sum_{k=1}^{N+2}b_k^2ke^{k-1}(1-(-1)^{k-1})=0,
\label{6_6}
\end{equation}
\begin{equation}
\sum_{k=0}^{N}a_k^1e^{k+1}\frac{(1-(-1)^{k+1})}{k+1}-\sum_{k=0}^{N}a_k^2e^{k+1}\frac{(1-(-1)^{k+1})}{k+1}=0,
\label{6_7}
\end{equation}
\begin{equation}
\sum_{k=0}^{N+2}b_k^1e^{k+1}\frac{(1-(-1)^{k+1})}{k+1}-\sum_{k=0}^{N+2}b_k^2e^{k+1}\frac{(1-(-1)^{k+1})}{k+1}=0.
\label{6_8}
\end{equation}

Hence, the equations (\ref{6_4})-(\ref{6_8}) form a system of $6N+10$ real linear algebraic equations with $6N+10$ real unknowns. Solving this system gives the Taylor polynomial approximations (\ref{6_1}) of the unknown functions.

The results obtained by this method have been verified on several examples by additionally solving the system (\ref{3_7}) numerically using spline approximations of the unknown functions. To obtain these approximations, first of all, rewrite the system (\ref{3_7}) in the form which involves integration over finite intervals only:
$$
\alpha_1\psi_1(x)+\beta_1\int_{|r|<l}\frac{\varphi_1(r)dr}{r-x}+\beta_1\int_{|r|<l}\frac{\varphi_0(r)dr}{rx-l^2}+\tau=\gamma_1^+\psi_1''(x)+f^+(x),\,\,\,|x|<l,
$$
$$
-\alpha_1\varphi_1(x)+\beta_1\int_{|r|<l}\frac{\psi_1(r)dr}{r-x}+\beta_1\int_{|r|<l}\frac{\psi_0(r)dr}{rx-l^2}+\sigma=-\gamma_0^+\psi_1'(x)+g^+(x),\,\,\,|x|<l,
$$
$$
\alpha_2\psi_2(x)-\beta_2\int_{|r|<l}\frac{\varphi_2(r)dr}{r-x}-\beta_2\int_{|r|<l}\frac{\varphi_0(r)dr}{rx-l^2}+\tau=-\gamma_1^-\psi_2''(x)+f^-(x),\,\,\,|x|<l,
$$
\begin{equation}
-\alpha_2\varphi_2(x)-\beta_2\int_{|r|<l}\frac{\psi_2(r)dr}{r-x}-\beta_2\int_{|r|<l}\frac{\psi_0(r)dr}{rx-l^2}+\sigma=\gamma_0^-\psi_2'(x)+g^-(x),\,\,\,|x|<l,
\label{6_9}
\end{equation}
$$
(\alpha_1-\alpha_2)\psi_0(x)+\beta_1l^2\int_{|r|<l}\frac{\varphi_1(r)dr}{rx-l^2}+\beta_2l^2\int_{|r|<l}\frac{\varphi_2(r)dr}{rx-l^2}+(\beta_1+\beta_2)\int_{|r|<l}\frac{\varphi_0(r)dr}{r-x}
$$
$$
=\gamma_1^il^{-4}(2x^2\psi_0(x)+4x^3\psi'_0(x)+x^4\psi''_0(x)),\,\,\,|x|<l,
$$
$$
-(\alpha_1-\alpha_2)\varphi_0(x)+\beta_1l^2\int_{|r|<l}\frac{\psi_1(r)dr}{rx-l^2}+\beta_2l^2\int_{|r|<l}\frac{\psi_2(r)dr}{rx-l^2}+(\beta_1+\beta_2)\int_{|r|<l}\frac{\psi_0(r)dr}{r-x}
$$
$$
=\gamma_0^il^{-2}(x\psi_0(x)+x^2\psi'_0(x)),\,\,\,|x|<l,
$$
where $\varphi_0(x)=\frac{l^2}{x}\varphi\left(\frac{l^2}{x} \right)$, $\psi_0(x)=\frac{l^2}{x}\psi\left(\frac{l^2}{x} \right)$.

Take the linear spline approximations of the functions $\varphi_1(x)$, $\varphi_2(x)$ and $\varphi_0(x)$:
\begin{equation}
\varphi_k(x)=\frac{w_k^j}{h}(x-x_{j-1})+\frac{w_k^{j-1}}{h}(x_j-x), \,\,\,x\in[x_{j-1},x_j],\,\,\,j=1,\ldots,2N,\,\,\,k=0,1,2,
\label{6_10}
\end{equation}
where the points $x_j=-l+jl/N$, $j=0,\ldots,2N$, generate the partition of the segment $[-l,l]$, and  the values of the functions $\varphi_k(x_j)$ at these points are equal to $w_k^j$. To guarantee the proper decline of the function $\varphi(x)$ at infinity (and, hence, the existence of all the corresponding integrals) it is necessary that $w_0^N=0$. Thus, the approximations (\ref{6_10}) contain $6N+2$ unknowns $w_k^j$.

To be able to differentiate the approximations of the functions $\psi_1(x)$, $\psi_2(x)$ and $\psi_0(x)$ twice it is necessary to take the higher order approximations of these functions which can be achieved by approximating the functions $\psi_1(x)$, $\psi_2(x)$ and $\psi_0(x)$ by cubic splines \cite{Schmidt1987}:
\begin{equation}
\psi_k(x)=S_k^j(x),\,\,\,x\in[x_{j-1},x_j],\,\,\,j=1,\ldots 2N,
\label{6_11}
\end{equation}
\begin{equation}
S_k^j(x)=\frac{z_k^j}{6h}(x-x_{j-1})^3+\frac{z_k^{j-1}}{6h}(x_j-x)^3
\label{6_12}
\end{equation}
$$
+\left(\frac{y_k^j}{h}-\frac{h}{6}z_k^j \right)(x-x_{j-1})+\left(\frac{y_k^{j-1}}{h}-\frac{h}{6}z_k^{j-1} \right)(x_j-x).
$$
It is easy to see that the functions $\psi_k(x)$ represented by (\ref{6_11}) and their second derivatives $\psi''_k(x)$ are continuous through the points $x_j$, $j=1,\ldots,(2N-1)$. To achieve the continuity of the derivatives $\psi'_k(x)$ at these points we must satisfy the following conditions
\begin{equation}
z_k^{j-1}+4z_k^j+z_k^{j+1}=\frac{6}{h^2}(y_k^{j-1}-2y_k^j+y_k^{j+1}),\,\,\,j=0,\ldots,2N.
\label{6_13}
\end{equation}
Also assume that in the formulas (\ref{6_12}), (\ref{6_13})
$$
z_k^{-1}=z_k^{2N+1}=0,\,\,\,y_k^{-1}=y_k^{2N+1}=0,\,\,\,k=0,1,2.
$$
Similarly, to guarantee an adequate decline of the function $\psi(x)$ at infinity it is necessary to set $y_0^N=0$. Thus, the formulas (\ref{6_11}), (\ref{6_12}) contain $12N+5$ unknowns $y_k^j$, $z_k^j$ and  the equations (\ref{6_13}) provide $6N+3$ conditions to find these unknowns.

Substituting the representations (\ref{6_10}), (\ref{6_11}) into the system (\ref{6_9}) and into the additional conditions (\ref{4_1})-(\ref{4_4}) produces the remaining $12N+4$ linear equations. Observe, that the equations of the system (\ref{6_9}) are satisfied numerically only at the $2N$ points $\xi_j=(x_{j-1}+x_j)/2$, $j=1,\ldots, 2N$, which are the midpoints of the segments $[x_{j-1},x_j]$. Solving the resulting system of linear algebraic equations numerically results in the spline approximations of the functions $\varphi_k(x)$, $\psi_k(x)$ at the points $x_j$, $j=0,\ldots,2N$.

\begin{figure}
	\centering
		\scalebox{0.5}{\includegraphics{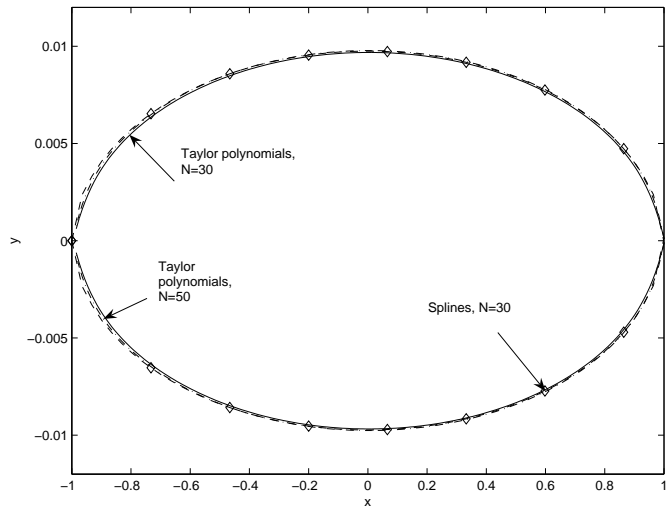}}
	\caption{Comparison of the results obtained by the Taylor polynomials and the spline approximations.}
	\label{fig2}
\end{figure}

Comparison of the results obtained by the approximations with the Taylor polynomials and with splines is presented on the fig. \ref{fig2}. The results have been obtained for $\mu_1=\mu_2=70$, $\nu_1=\nu_2=0.3$, $l=1$, and the following surface tension constants: $\gamma_1^+=\gamma_1^-=0.01$, $\gamma_0^+=\gamma_0^-=0.01$, $\gamma_1^i=\gamma_0^i=0.005$. The construction is subjected to the pure tensile loading $\sigma=1$ at infinity and the boundaries of the crack are free from external stresses. The solid line depicts the shape of the crack obtained with Taylor polynomial approximations for $N=30$, the dotted line represents the results obtained with Taylor polynomial approximations for $N=50$, and, finally, the diamond marks represent the spline approximations for $N=30$. It is easy to see that all the graphs are almost indistinguishable from each other. Since the Taylor polynomial approximations were, in author's opinion, more computationally efficient, the results presented further in this paper were obtained by this method. 

\begin{figure}
	\centering
		\scalebox{0.5}{\includegraphics{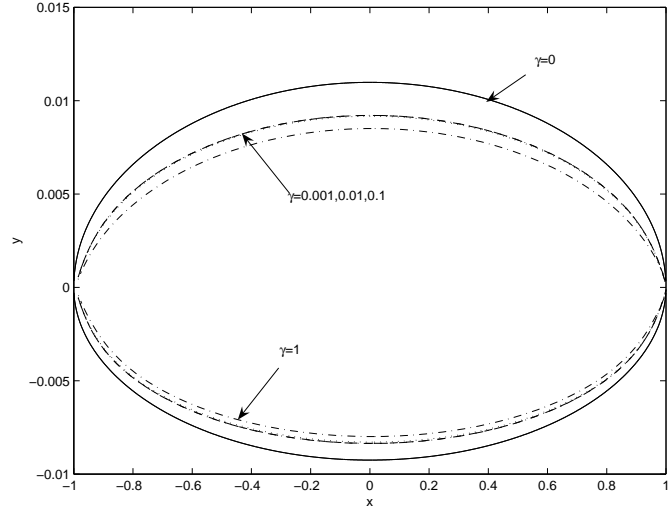}}
	\caption{Influence of the surface tension on the shape of the crack.}
	\label{fig3}
\end{figure}

The results for the interface crack problem in the case of the LEFM are well-known and have been derived by many authors for different types of configurations and loading. In this paper, we compare the results obtained with the curvature-dependent surface tension model with the classical results of England \cite{England1965}. In the case of the absent loading at infinity and a uniform pressure $f^+(x)=f^-(x)=0$, $g^+(x)=g^-(x)=T$,  applied to the boundaries of the crack, the displacements of the boundaries of the crack are derived in an explicit form according to the following formulas \cite{England1965}:
\begin{equation}
u_2^+(x)=\frac{T(1+\kappa_1)\sqrt{\alpha}}{2\mu_1(1+\alpha)}(l^2-x^2)^{1/2}\cos\left(\gamma\ln\left|\frac{l+x}{l-x} \right| \right),\,\,\,|x|<l,
\label{6_14}
\end{equation}
\begin{equation}
u_2^-(x)=-\frac{T(1+\kappa_2)\sqrt{\alpha}}{2\mu_2(1+\alpha)}(l^2-x^2)^{1/2}\cos\left(\gamma\ln\left|\frac{l+x}{l-x} \right| \right),\,\,\,|x|<l,
\label{6_15}
\end{equation}
where
$$
\alpha=\frac{\mu_1+\mu_2\kappa_1}{\mu_2+\mu_1\kappa_1},\,\,\,\gamma=\frac{1}{2\pi}\ln\alpha.
$$
Observe the presence of the oscillating singularities given by the terms $\cos\left(\gamma\ln\left|\frac{l+x}{l-x} \right| \right)$ in the formulas (\ref{6_14}), (\ref{6_15}). 

\begin{figure}
	\centering
		\scalebox{0.72}{\includegraphics{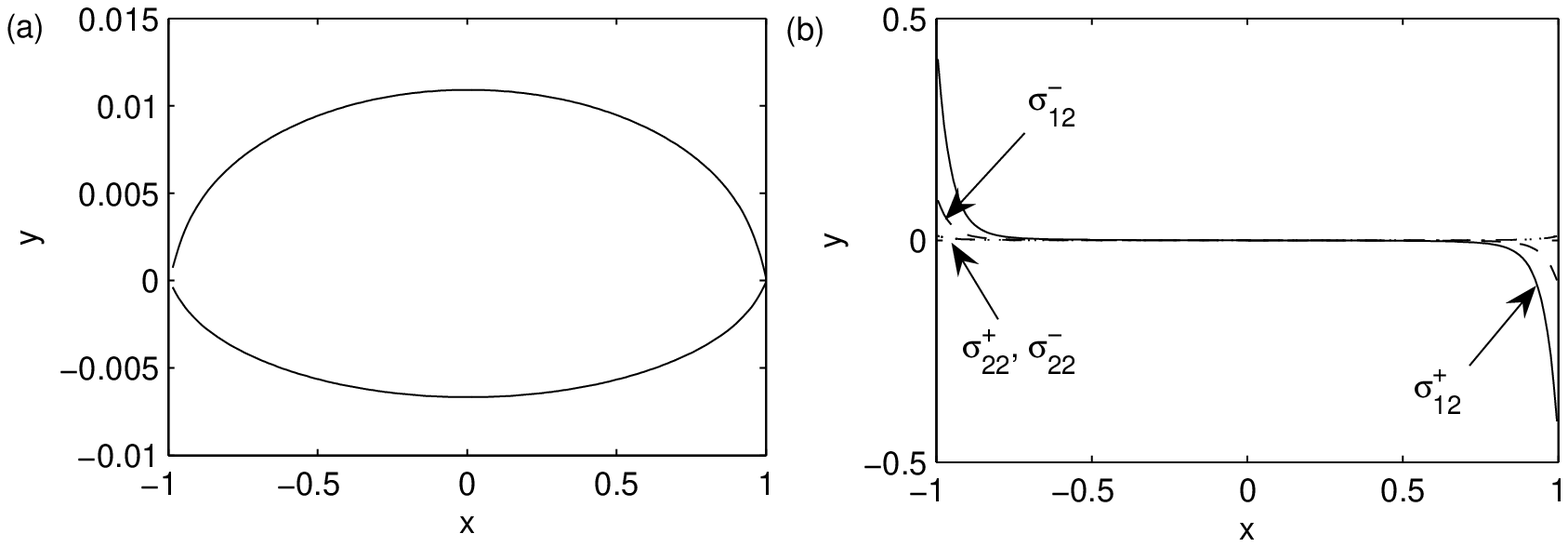}}
		\scalebox{0.7}{\includegraphics{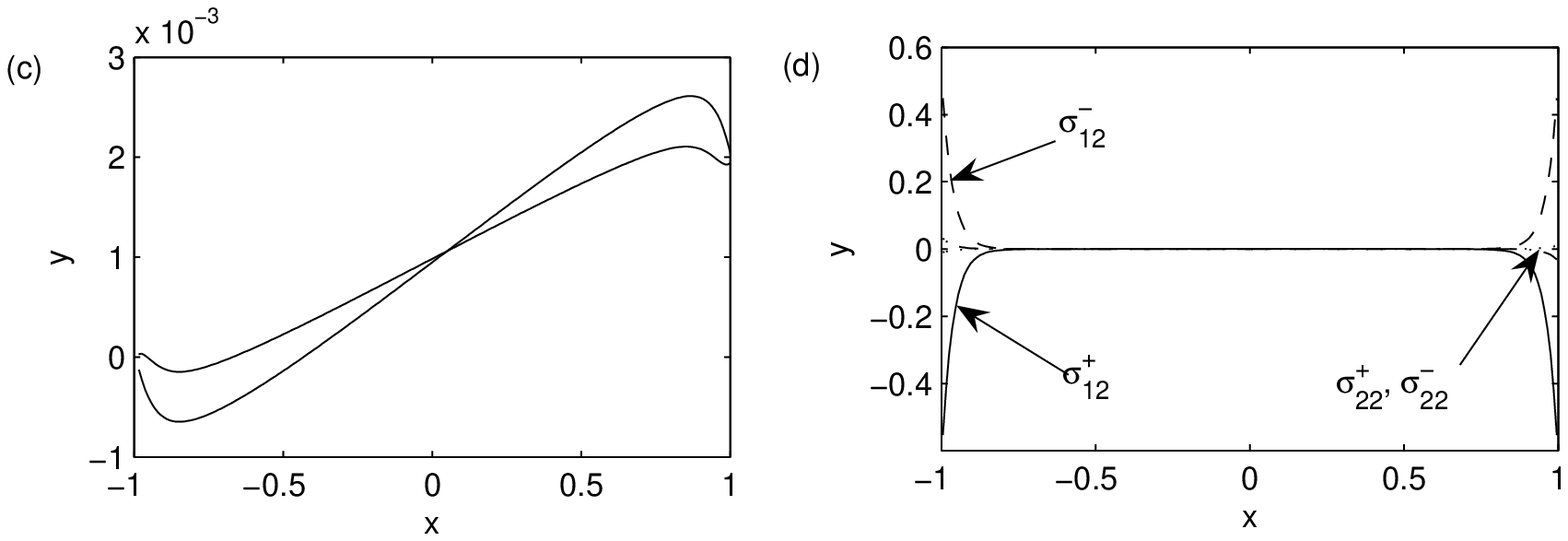}}
		\scalebox{0.7}{\includegraphics{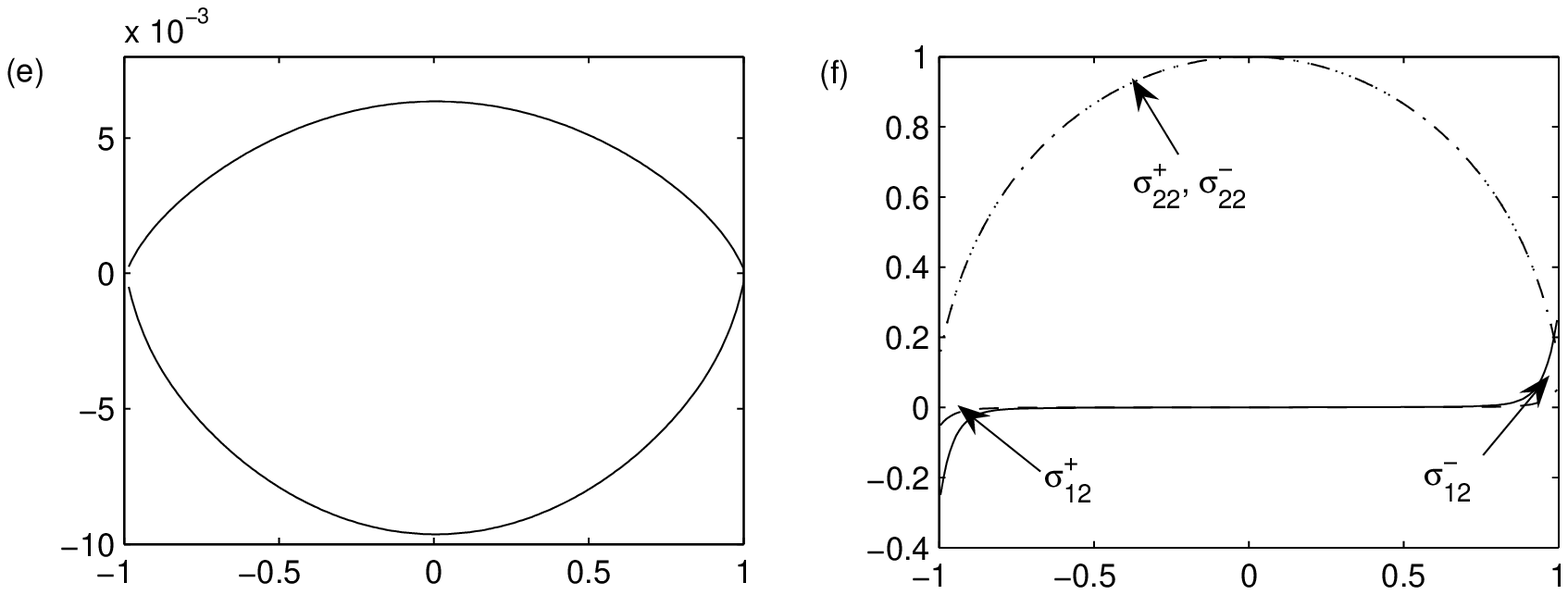}}				
	\caption{Shape of the crack and stresses on the crack boundary for different types of loading.}
	\label{fig4}
\end{figure}

The influence of the magnitude of the surface tension on the shape of the crack is presented on the fig. \ref{fig3}. The shape of the crack is plotted for $l=1$, $\mu_1=70$, $\mu_2=80$, $\nu_1=0.3$, $\nu_2=0.35$, and the following surface tension constants: $\gamma_1^+=\gamma_1^-=\gamma_1^i=\gamma_0^+=\gamma_0^-=\gamma_0^i=\gamma$. The constant $\gamma$ takes on the following five values: $\gamma=0$ (corresponds to the classical LEFM solution obtained by the formulas (\ref{6_14}), (\ref{6_15})), $\gamma=0.001$, $\gamma=0.01$, $\gamma=0.1$ and $\gamma=1.0$. Observe, that the introduction of the surface tension principally changes the behavior of the stresses and the displacements at the crack tips, however, the magnitude of the surface tension constants has little influence on the displacements of the boundaries of the crack. Observe also, that the displacements are almost equal for $\gamma=0.001$, $\gamma=0.01$ and $\gamma=0.1$, and only for $\gamma=1.0$ a noticeable difference can be seen. 

\begin{figure}[htbp]
	\centering
		\scalebox{0.7}{\includegraphics{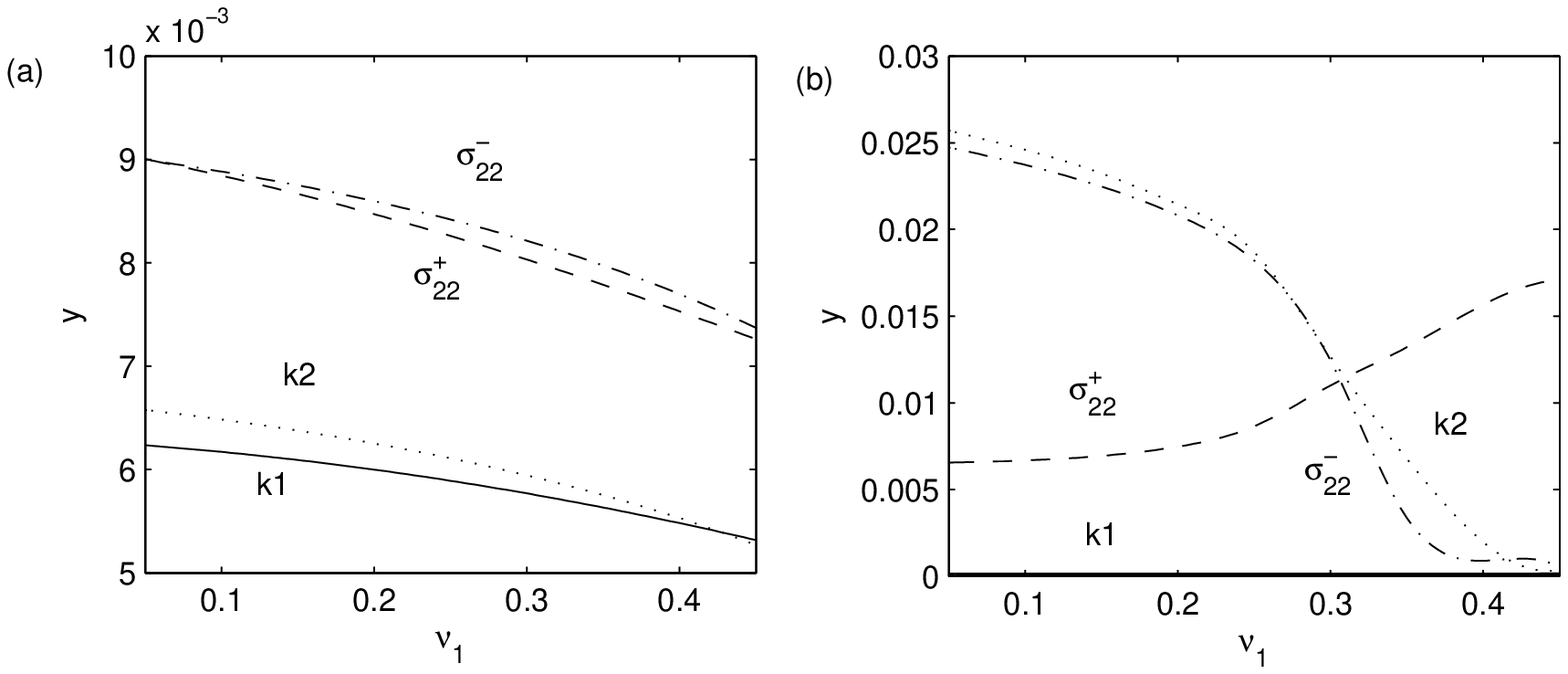}}	
	\caption{Dependence of the maximal stresses and the singularity coefficients on the Poisson's ratio $\nu_1$.}
	\label{fig5}
\end{figure}
\begin{figure}[htbp]
	\centering
		\scalebox{0.7}{\includegraphics{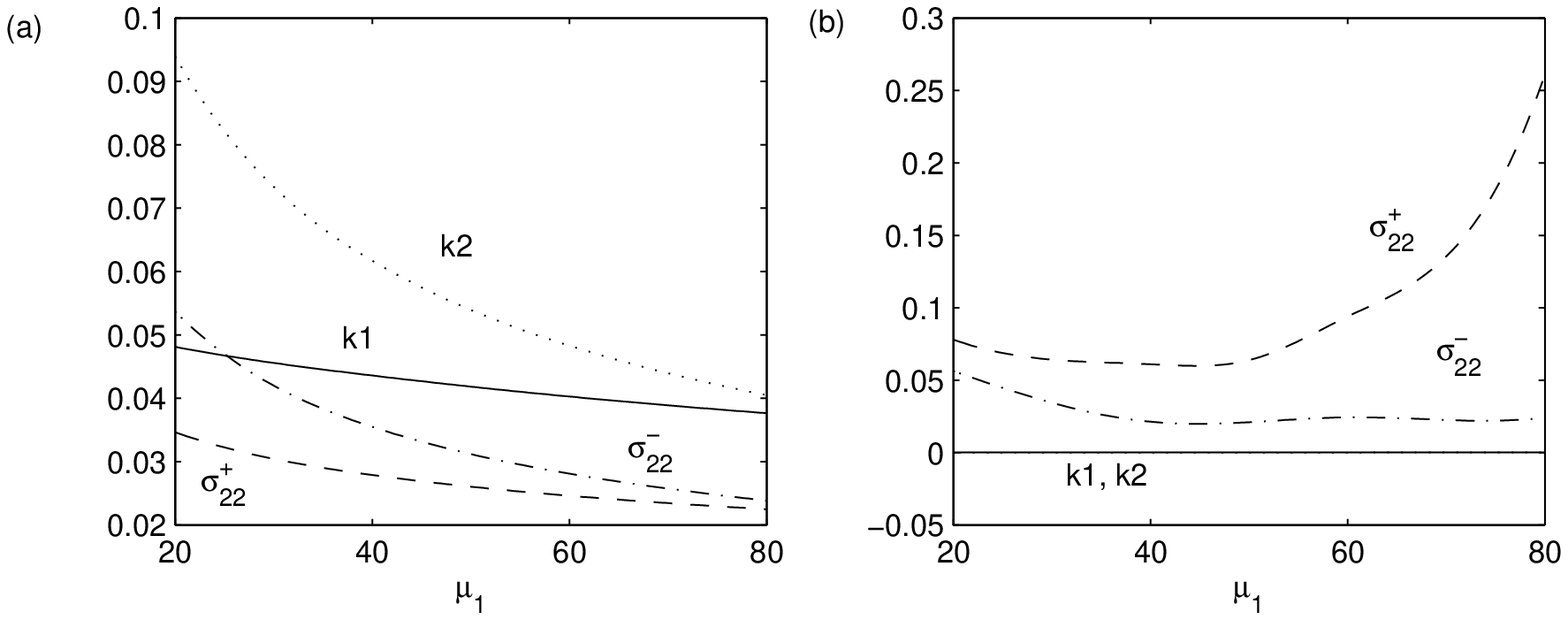}}	
	\caption{Dependence of the maximal stresses and the singularity coefficients on the shear modulus $\mu_1$.}
	\label{fig6}
\end{figure}

The shape of the crack and the stresses on the crack boundary are shown on the fig. \ref{fig4} for different types of loading of the construction. The computations are made for $\mu_1=70$, $\mu_2=80$, $\nu_1=0.3$, $\nu_2=0.35$, $l=1$, $\gamma_1^+=0.02$, $\gamma_1^-=0.01$, $\gamma_0^+=0.01$, $\gamma_0^-=0.02$, $\gamma_1^i=0.01$, $\gamma_0^i=-0.01$. The graphs (a), (b) correspond to the shape of the crack and the stresses on the crack boundary for the case when the construction is subjected to the pure tensile loading $\sigma=1$ at infinity in the absence of any other loading; the graphs (c), (d) to the case of the pure shear loading at infinity $\tau=1$; the graphs (e), (f) to the case of the loading $f^+(x)=f^-(x)=0$, $g^+(x)=g^-(x)=\sqrt{l^2-x^2}$ applied on the boundary of the crack and the absence of the loading at infinity. Note that the edges of the crack on the graph (c) self-intersect. This is because of the appearance of the contact zones on the boundaries of the crack in the case of the shear loading. The self-intersection may be removed if these zones are explicitly accounted for during the solution of the problem. Also observe that with the exception of the stresses $\sigma_{22}^{\pm}$ on the fig. \ref{fig4}(f), the stresses on the boundaries of the crack are close to zero everywhere except for the small area near the tips of the crack. This is due to the fact that the curvature dependent surface tension is small in magnitude everywhere except for the small zones near the crack tips. However, presence of this tension essentially changes the character of the singularities at the crack tips.

\begin{figure}
	\centering
		\scalebox{0.7}{\includegraphics{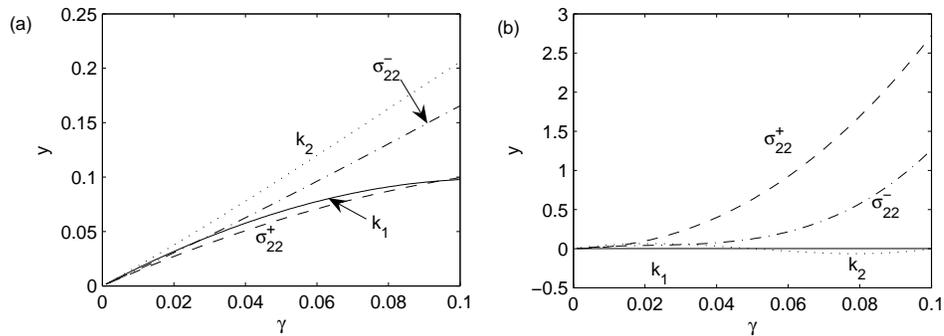}}	
	\caption{Dependence of the maximal stresses and the singularity coefficients on the parameter $\gamma$.}
	\label{fig7}
\end{figure}

The dependence of the maximal stresses and the singularity coefficients in the construction of the Poisson's ratio of one of the semi-planes $\nu_1$ is shown of the fig. \ref{fig5}. The other parameters of the construction are $\mu_1=\mu_2=70$, $\nu_2=0.30$, $l=1$, $\gamma_1^+=\gamma_1^-=\gamma_0^+=\gamma_0^-=0.01$,  $\gamma_1^i=\gamma_0^i=0.005$. The figure (a) corresponds to the pure tensile loading at infinity $\sigma=1$ and the figure (b) corresponds to the pure shear loading $\tau=1$. The coefficients $k_1$ and $k_2$ are the coefficients by the square logarithmic singularities of the stresses $\sigma_{12}^{\pm}$:
$$
\sigma_{12}^+(x)=k_1\ln^2(l-x)+O(\ln(l-x)),\,\,\, \sigma_{12}^-(x)=k_2\ln^2(l-x)+O(\ln(l-x)).
$$
Observe that the maximums of the stresses $\sigma_{22}^{\pm}$ are achieved at the tips of the crack. Similar graphs for the dependence of the maximal stresses and the singularity coefficients on the shear modulus $\mu_1$ are shown on the fig. \ref{fig6}. All the parameters are the same as before with the exception of the changing parameter $\mu_1$, and $\nu_1=0.3$. The graphs of the dependence of the maximal stresses and the singularity coefficients on the parameter $\gamma$, $\gamma_1^+=\gamma_1^-=\gamma_0^+=\gamma_0^-=\gamma$, $\gamma_1^i=\gamma_0^i=0.005$ are shown on the fig. \ref{fig7}. The other parameters of the construction are $\mu_1=70$, $\mu_2=80$, $\nu_1=0.3$, $\nu_2=0.35$, $l=1$. The figure (a) as before correspond to the pure tensile loading $\sigma=1$ at infinity and the figure (b) to the pure shear loading $\tau=1$.

\section*{Conclusions}

In this paper, a new model of brittle fracture first proposed in \cite{Slatteryetal} and later used to study fracture problems in \cite{OhWaltonSlattery2006}-\cite{SendovaWalton2010a} is applied to the problem of an interface fracture between two dissimilar isotropic semi-planes under the tensile-shear in-plane loading. The classical boundary condition of linear elasticity is augmented with a curvature-dependent surface tension. 

The physical problem is reduced to the system of singular integro-differential equations by applying Dirichlet-to-Neumann mappings. The system of singular integro-differential equations is further reduced to the system of Fredholm integral equations. It is shown that unlike in the classical solutions of linear elasticity, the introduction of the surface tension leads to the solutions which do not possess integrable power singularities of the order $1/2$ or the oscillating singularities, thus making the solutions more practically feasible. The integrable singularities of the logarithmic type may still be present. The most important problem arising from the presented solution is in developing more general models of surface tension which may allow to eliminate logarithmic singularities as well, thus making the solutions completely bounded in the neighborhoods of the crack tips. This is the subject of ongoing research.

The numerical solution of the system is obtained by approximating the unknown functions by Taylor polynomials and is further verified by spline collocation method. The results of computations are presented for different types of mechanical parameters.

\vspace{.1in}

\section*{Acknowledgments} 
The author is grateful to Prof. J. R. Walton for the suggestion of the topic and many helpful discussions.

This publication is based on work supported by Award No. KUS-C1-016-04, made by King Abdullah University of Science and Technology (KAUST).

\end{document}